\documentclass[11pt]{article}
\usepackage{epsfig}

\begin{document}

\title{
Discrete Riccati equation, hypergeometric functions and circle
patterns of Schramm type}

\author{{\Large
Agafonov S.I. } \\
\\
    Department of Mathematical Sciences \\
Loughborough University \\
    Loughborough, Leicestershire LE11 3TU \\
United Kingdom \\
    e-mail: {\tt
SAgafonov@rusfund.ru}
}
\date{}
\maketitle

\unitlength=1mm

\newtheorem{theorem}{Theorem}
\newtheorem{proposition}{Proposition}
\newtheorem{lemma}{Lemma}
\newtheorem{corollary}{Corollary}
\newtheorem{definition}{Definition}

\pagestyle{plain}

\begin{abstract}

\noindent Square grid circle patterns with  prescribed
intersection angles, mimicking holomorphic
 maps $z^{\gamma }$ and ${\rm log}(z)$  are studied.
It is shown that the corresponding
 circle patterns are embedded and described by
special separatrix solutions of  discrete  Painlev\'e and Riccati
equations. General solution of this Riccati equation is expressed
in terms of the hypergeometric function. Global properties of
these solutions, as well as of the discrete $z^{\gamma }$ and
${\rm log}(z)$, are established.
\end{abstract}

\section{Introduction}
The theory of circle patterns is a rich fascinating area having
its origin in classical theory of circle packings. Its fast
development in recent years is caused by mutual influence and
interplay of ideas and concepts from discrete geometry, complex
analysis and the theory of integrable systems.

The progress in this area was initiated by Thurston's idea
\cite{T},\cite{MR} about approximating the Riemann mapping by
circle packings. Classical circle packings comprised of disjoint
open disks were later generalized to circle patterns where the
disks may overlap (see for example \cite{H}). Different underlying
combinatorics were considered. Schramm introduced a class of
circle patterns with the combinatorics of the square grid
\cite{Schramm}; hexagonal circle patterns were studied in
\cite{BHS} and \cite{BH}.

The striking analogy between circle patterns and the classical
analytic function theory is underlined by such facts as the
uniformization theorem
 concerning circle packing realizations of
cell complexes of a prescribed combinatorics \cite{BS}, discrete
maximum principle, Schwarz's lemma \cite{R} and rigidity
properties \cite{MR},\cite{H}, discrete Dirichlet principle
\cite{Schramm}.

The convergence of discrete conformal maps represented by circle
packings was proven by Rodin and Sullivan \cite{RS}. For a
prescribed regular combinatorics this result was refined. He and
Schramm \cite{HS} showed that  for hexagonal packings the
convergence is $C^{\infty}.$ The uniform convergence for circle
patterns with the combinatorics of the square grid and orthogonal
neighboring circles was established by Schramm \cite{Schramm}.

Approximation issue naturally leads to the question about analogs
to standard holomorphic functions. Computer experiments give
evidence for their existence \cite{DS},\cite{TH} however not very
much is known.  For circle packings with
 the hexagonal combinatorics the only explicitly described examples are Doyle spirals  \cite{Doy},\cite{BDS} which are  discrete
analogs of exponential maps and conformally symmetric packings,
which are analogs of a quotient of Airy functions \cite{BHConf}.
For patterns with overlapping circles more explicit examples are
known: discrete versions of ${\rm exp} (z)$, ${\rm erf}(z)$
\cite{Schramm},  $z^{\gamma}$, $ {\rm log} (z)$ \cite{AB} are
constructed for patterns with underlying combinatorics of the
square grid; $z^{\gamma}$, ${\rm log}(z)$ are also described for
hexagonal patterns \cite{BHS}, \cite{BH}.

It turned out that an effective approach to the description  of
circle patterns with overlapping circles is given by the theory of
integrable systems  (see \cite{BPD},\cite{BHS},\cite{BH}). For
example, Schramm's circle patterns are governed by a difference
equation which is the stationary Hirota equation \cite{Schramm}
(see \cite{Z} for a survey). This approach proved to be especially
useful for the construction of discrete $z^{\gamma}$ and ${\rm
log}(z)$ in \cite{AB},\cite{BHS},\cite{BH} with the aid of some
isomonodromy problem. Another connection with the theory of
discrete integrable equations was revealed in \cite{A},\cite{AB}:
embedded circle patterns are described by special solutions of
discrete Painlev\'e II equations, thus giving geometrical
interpretation thereof.

This research was motivated by the attempt to carry the results of
\cite{A},\cite{AB} over square grid circle patterns with
prescribed intersection angles giving Schramm patterns as a
special case. Namely, we prove that such circle patterns mimicking
$z^{\gamma}$ and ${\rm log}(z)$ are embedded. This turned out to
be not straightforward and lead to asymptotical analysis of
solutions to {\it discrete Riccati equation}. As the Riccati
differential equation is known as the only equation of the first
order possessing Panlev\'e property we are tempted to conclude
that circle patterns are described by discrete equations with
Painlev\'e property though there is no satisfactory generalization
thereof on discrete equations.

\smallskip

\begin{figure}[th]
\begin{center}
\epsfig{file=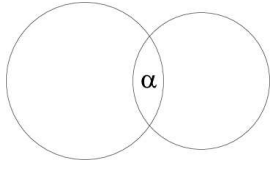,width=35mm}
 \caption{Circles intersection angle.}
\label{Angle}
\end{center}
\end{figure}

\noindent We use the following definition for square grid circle
patterns, which is slightly modified version of one from
\cite{Schramm}.
\begin{definition}
Let $G$ be a subgraph of the 1-skeleton of the cell complex with
vertices ${\bf Z}+i{\bf Z}={\bf Z^2} $ and $0< \alpha <\pi $.
Square grid circle pattern for $G$ with intersection angles
$\alpha$ is an indexed collection of circles on the complex plane
$$ \{C_z:z\in V(G),\ C_z\in {\bf C} \} $$ that satisfy: \\ 1) if
$z,z+i \in V(G)$ then the intersection angle of $C_z$ and
$C_{z+i}$ is $\alpha$ ,\\ 2) if $z,z+1 \in V(G)$ then the
intersection angle of $C_z$ and $C_{z+1}$ is $\pi-\alpha$,\\
\smallskip
3) if $z,z+1+i \in V(G)$ (or $z,z-1+i \in V(G)$) then the disks,
defined by $C_z$ and $C_{z+1+i}$ ($C_z$ and $C_{z-1+i}$
respectively)
 are tangent and disjoint,\\
4) if $z,z_1,z_2\in V(G)$, $|z_1-z_2|=\sqrt{2}$,
$|z-z_1|=|z-z_2|=1$ (i.e. $C_{z_1}$,$C_{z_2}$ are tangent and
$C_z$ intersects $C_{z_1}$ and $C_{z_2}$) and $z_2=z+i(z_1-z)$
(i.e. $z_2$ is one step counterclockwise from $z_1$), then the
circular order of the triplet of points $C_z\cap
C_{z_1}-C_{z_2}$,$C_{z_1}\cap C_{z_2}$,$C_z\cap C_{z_2}-C_{z_1}$
agrees with the orientation of $C_z$.
\end{definition}
The intersection angle is the angle at the corner of disc
intersection domain (Fig. \ref{Angle}).

\begin{figure}[th]
\begin{center}
\epsfig{file=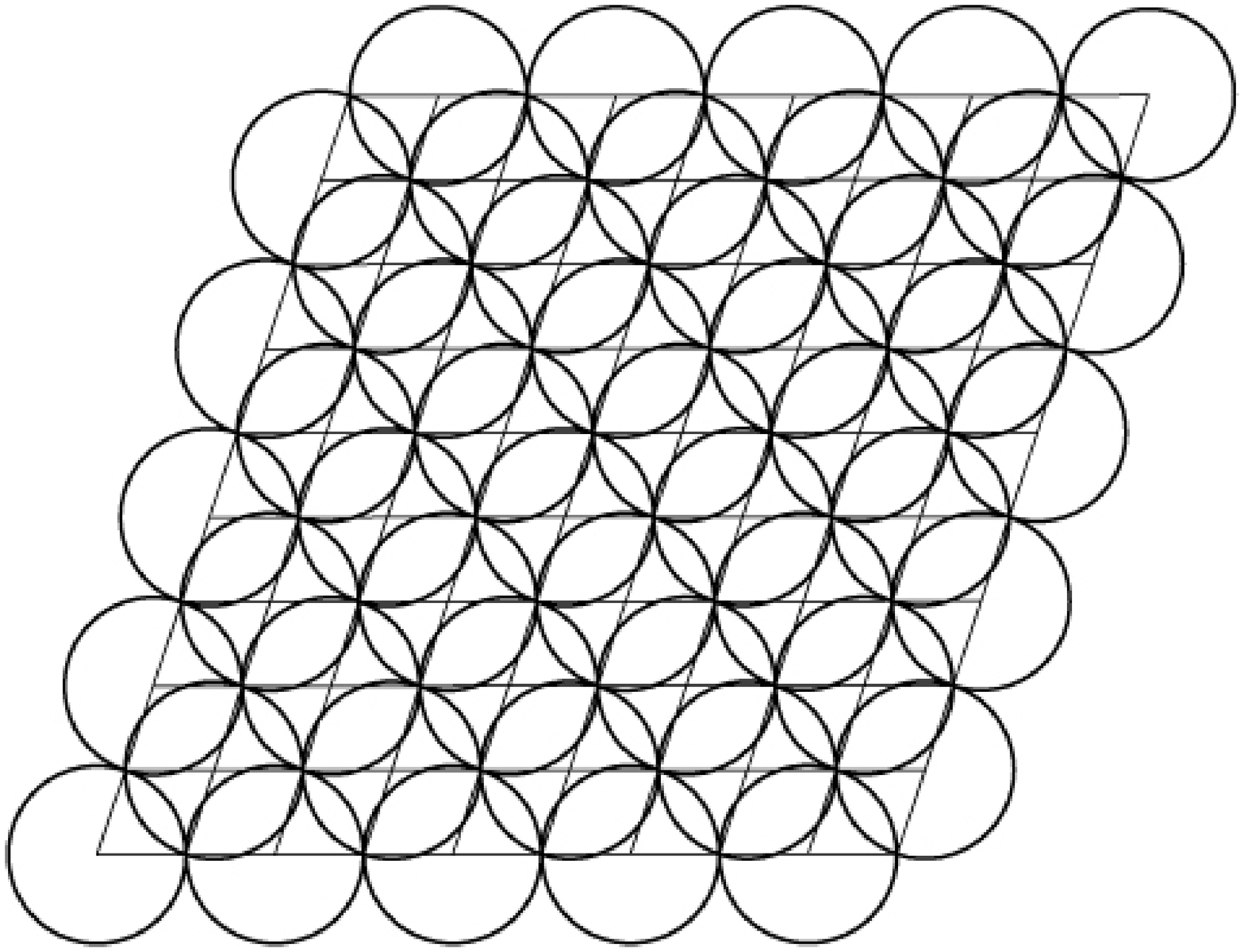,width=45mm}
\begin{picture}(20,50)
\put(2,22){\vector(1,0){15}} \put(7,27){ \it  \huge  f}
\end{picture}
\epsfig{file=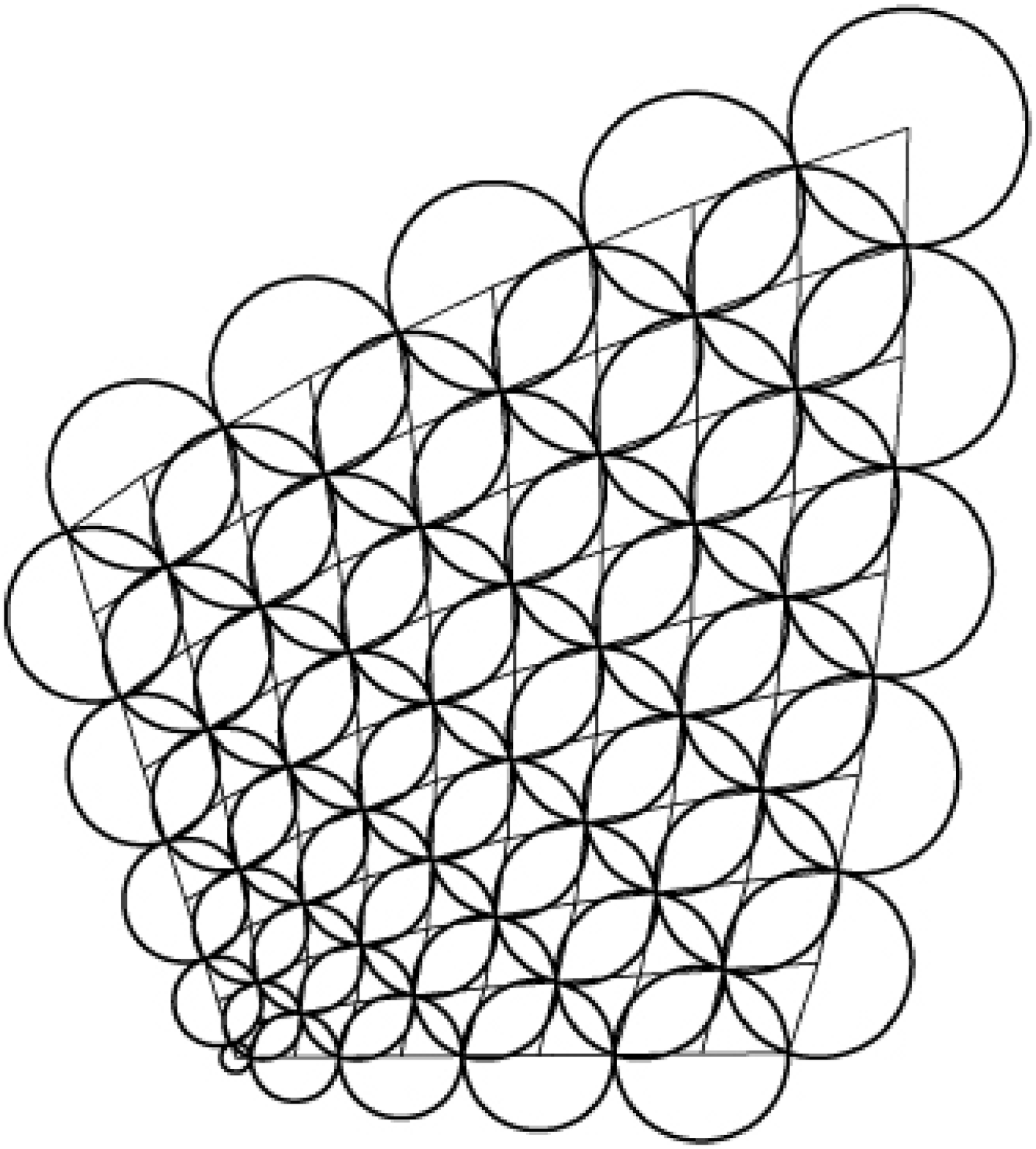,width=50mm} \caption{Schramm type circle
patterns as a discrete conformal map.  The discrete version of the
holomorphic mapping $z^{3/2}$.  The case $\tan \alpha =3.$}
\label{CPMap}
\end{center}
\end{figure}

\noindent  To visualize the analogy between Schramm's circle
patterns and conformal maps, consider regular patterns composed of
unit circles and suppose that the radii are being deformed to
preserve the intersection angles of neighboring circles and the
tangency of half-neighboring ones. Discrete maps taking the
intersection points and the centers of the unit circles of the
standard regular patterns to the respective points of the deformed
patterns mimic classical holomorphic functions, the deformed radii
being analog of $|f'(z)|$ (see Fig. \ref{CPMap}).

In section \ref{CP} we give definition of discrete $Z^{\gamma}$ as
a solution to some integrable equation subjected to a
non-autonomous constraint. Its geometrical properties of immersion
and embeddedness are expressed in terms of solution for radii of
corresponding circle patterns.

For this solution to be positive it is necessary that some
discrete Riccati equation has positive solution. This equation is
studied in section \ref{RE} where its general solution is
expressed via  the hypergeometric function.

Section \ref{IPE} completes the proof of embeddednes,  discrete
equations of Painlev\'e type being the main tool. Possible
generalizations for non-regular combinatorics and for non-circular
patterns are discussed in section \ref{CR}.

\section{ Discrete $Z^{\gamma}$ and square grid circle patterns of Schramm
type.} \label{CP}
\begin{definition}\label{skew-Schramm-def}
Discrete map $Z^{\gamma}$, $0<\gamma <2$ is the solution $f:{\bf
Z^2_+} \to {\bf C}$ of
\begin{equation}\label{q}
q(f_{n,m},f_{n+1,m},f_{n+1,m+1},f_{n,m+1})=e^{-2i\alpha }
\end{equation}
\begin{equation}\label{sg-c}
\gamma
f_{n,m}=2n\frac{(f_{n+1,m}-f_{n,m})(f_{n,m}-f_{n-1,m})}{f_{n+1,m}-f_{n-1,m}}+
2m\frac{(f_{n,m+1}-f_{n,m})(f_{n,m}-f_{n,m-1})}{f_{n,m+1}-f_{n,m-1}},
\end{equation}
with $0<\alpha<\pi$ and the initial conditions
\begin{equation} \label{right-sg-initial}
  f_{1,0}=1,\
  f_{0,1}=e^{i \gamma \alpha  },
\end{equation}
where $q$ stands for cross-ratio of elementary quadrilaterals:
$$q(f_1,f_2,f_3,f_4)=\frac{(f_1-f_2)(f_3-f_4)}
{(f_2-f_3)(f_4-f_1)}.$$
\end{definition}
Constraint (\ref{sg-c}) was obtained from some isomonodromy
problem of Lax representation of (\ref{q}), found for $\gamma=1$
in  \cite{NC} (see also \cite{AB},\cite{BH}). The definition of
$Z^{\gamma}$ can be justified by the following properties:
\begin{itemize}
\item if one thinks of $f$ as defined on the vertices of the cell complex with
diamond-shaped faces (see Fig.\ref{CPMap}) then (\ref{q}) means
that $f$ respects the cross-ratios of the faces and therefore is
"locally conformal",
\item asymptotics of (\ref{sg-c}) as $n,m \to \infty$ suggests that $f$ approximates
$z^{\gamma}$.
\end{itemize}

\medskip

\noindent {\bf Remark.} Equation (\ref{q}) with $\alpha =\pi/2 $
was used in \cite{BPdis} to define {\it discrete conformal maps}.
The motivation was that $f$ maps vertices of squares into vertices
of the "conformal squares". Consider the surface glued of these
conformal squares along the corresponding edges. This surface is
locally flat but can have cone-like singularities in vertices. If
the map is an immersion then the corresponding surface does not
have such singularities. Therefore it is more consistent to define
as discrete conformal an {\it immersion} map on the vertices of
cell decomposition of ${\it C}$ which preserves cross-ratios of
its faces.

\begin{proposition}\label{compatible}\cite{BH, AB}
Constraint (\ref{sg-c}) is compatible with (\ref{q}).
\end{proposition}

Compatibility is understood as a solvability of some Cauchy
problem. In particular a solution to (\ref{q}),(\ref{sg-c}) in the
subset ${\bf Z^2_+}$ is uniquely determined by its values $
 f_{1,0},\ f_{0,1}$.  Indeed, constraint (\ref{sg-c}) gives $f_{0,0}=0$
and  defines $f$ along the coordinate axis $(n,0),\ (0,m)$ as a
second-order difference equation. Then all other $f_{k,m}$ with
$(k,m)\in {\bf Z^2_+}$ are calculated via  cross-ratios (\ref{q}).

In this paper the following initial conditions for
(\ref{q}),(\ref{sg-c}) are considered:
\begin{equation}\label{skew-Shramm-initial} f_{1,0}=1,\
 f_{0,1}=e^{i \beta      }
\end{equation} with real $\beta$. The solution $f$ defines a  circle pattern with square grid
combinatorics. For $\alpha =\pi /2$ such circle patterns were
introduced by Schramm in \cite{Schramm}. For any $0<\alpha <\pi$
we obtain a natural generalization of Schramm circle patterns.

In what follows we say that the triangle $(z_1,z_2,z_3)$ has {\it
positive (negative) orientation} if
$$\frac{z_3-z_1}{z_2-z_1}=\left|\frac{z_3-z_1}{z_2-z_1}\right|e^{i\phi}
\ \ {\rm with} \ 0\le  \phi \le \pi \  \ \ (-\pi <\phi < 0).$$

\begin{lemma}\label{geom-kite}
Let $q(z_1,z_2,z_3,z_4)=e^{-2i\alpha }$,  $0<\alpha <\pi$.
\begin{itemize}
\item If $|z_1-z_2|=|z_1-z_4|$ and the triangle $(z_1,z_2,z_4)$ has
positive orientation then $|z_3-z_2|=|z_3-z_4|$ and the angle
between $[z_1,z_2]$ and $[z_2,z_3]$ is $(\pi-\alpha )$.
\item If $|z_1-z_2|=|z_1-z_4|$ and the triangle $(z_1,z_2,z_4)$ has
negative orientation then $|z_3-z_2|=|z_3-z_4|$ and the angle
between $[z_1,z_2]$ and $[z_2,z_3]$ is $\alpha$.
\item If the angle
between $[z_1,z_2]$ and $[z_1,z_4]$ is $\alpha$ and the triangle
$(z_1,z_2,z_4)$ has positive orientation then
$|z_3-z_2|=|z_1-z_2|$ and   $|z_3-z_4|=|z_4-z_1|$.
\item If the angle
between $[z_1,z_2]$ and $[z_1,z_4]$ is $(\pi- \alpha)$ and the
triangle $(z_1,z_2,z_4)$ has negative orientation then
$|z_3-z_2|=|z_1-z_2|$ and   $|z_3-z_4|=|z_4-z_1|$.
\end{itemize}
\end{lemma}
{\it Proof:} straightforward.

\begin{proposition} \label{kite}
All the elementary quadrilaterals $(f_{n,m},
f_{n+1,m},f_{n+1,m+1},f_{n,m+1})$ for the solution of
(\ref{q}),(\ref{sg-c}) with initial (\ref{skew-Shramm-initial})
are of kite form: all edges at the vertex $f_{n,m}$ with $n+m=0 \
({\rm mod} \ 2)$ are of the same length. Moreover, each elementary
quadrilateral has one of the forms enumerated in lemma
\ref{geom-kite}.
\end{proposition}
\noindent {\it Proof:}   Given initial $f_{0,1}$ and $f_{1,0}$
constraint (\ref{sg-c}) gives $f_{n,0}$ and $f_{0,m}$ for all $n,m
\ge 1.$ It is easy to check that $f$ has the following equidistant
property:
\begin{equation} \label{equidistant}
f_{2n,0}-f_{2n-1,0}=f_{2n+1,0}-f_{2n,0},\ \
f_{0,2m}-f_{0,2m-1}=f_{0,2m+1}-f_{0,2m}
\end{equation}
for any $n\ge1$, $m\ge1$. Lemma \ref{geom-kite} and
$|f_{1,0}-f_{0,0}|=|f_{0,1}-f_{0,0}|$ allows one to apply
induction in $n,m$, starting with $n=0,\ m=0$.

Proposition \ref{kite} implies that for  $n+m=0 \ ({\rm mod} \ 2)$
the points $f_{n\pm 1,m}$,$f_{n,m\pm 1}$ lie on the circle with
the center at $f_{n,m}$. For the most $\beta $ (namely for $\beta
\ne \alpha$) the behavior of thus obtained circle pattern is quite
irregular: inner parts of different elementary quadrilaterals
intersect.

\begin{definition}
A discrete  map $f_{n,m}$ is called an immersion if inner parts of
adjacent elementary quadrilaterals
$(f_{n,m},f_{n+1,m},f_{n+1,m+1},f_{n,m+1})$ are disjoint.
\end{definition}

Consider the sublattice $\{n,m: \ n+m=0 \ ({\rm mod} \ 2)\}$ and
denote by $\bf V$ its quadrant $$ {\bf V}=\{z=N+iM:\ N,M \in {\bf
Z^2}, M \ge |N| \}, $$ where $$ N=(n-m)/2, \ \ M=(n+m)/2. $$
 We use complex labels $z=N+iM$ for this sublattice. Denote by
$C(z)$ the circle of the radius

\begin{equation}
R_z=|f_{n,m}-f_{n+1,m}|=|f_{n,m}-f_{n,m+1}|=|f_{n,m}-f_{n-1,m}|=|f_{n,m}-f_{n,m-1}|
\label{rmap}
\end{equation}
with the center at $f_{N+M,M-N}=f_{n,m}.$

Let $\{C(z)\}, \ z\in {\bf V}$ be a square grid circle pattern
 on the complex plane. Define $f_{n,m}: {\bf
Z^2_+ \to C}$ as follows:\\ a) if $n+m=0 \ ({\rm mod} \ 2)$ then
$f_{n,m}$ is the center of $C(\frac{n-m}{2}+i\frac{n+m}{2}),$\\ b)
if $n+m=1 \ ({\rm mod} \ 2)$ then
$
f_{n,m}:=C(\frac{n-m-1}{2}+i\frac{n+m-1}{2}) \cap
C(\frac{n-m+1}{2}+i\frac{n+m+1}{2})=
C(\frac{n-m+1}{2}+i\frac{n+m-1}{2}) \cap
C(\frac{n-m-1}{2}+i\frac{n+m+1}{2}).
$
Since all elementary quadrilaterals $(f_{n,m},
f_{n+1,m},f_{n+1,m+1},f_{n,m+1})$ are of kite form equation
(\ref{q}) is satisfied automatically. In what follows the function
$f_{n,m},$ defined as above by a) and b) is called {\it a discrete
 map corresponding to the circle pattern $\{C(z)\}$ .}

\begin{figure}[th]
\begin{center}
\epsfig{file=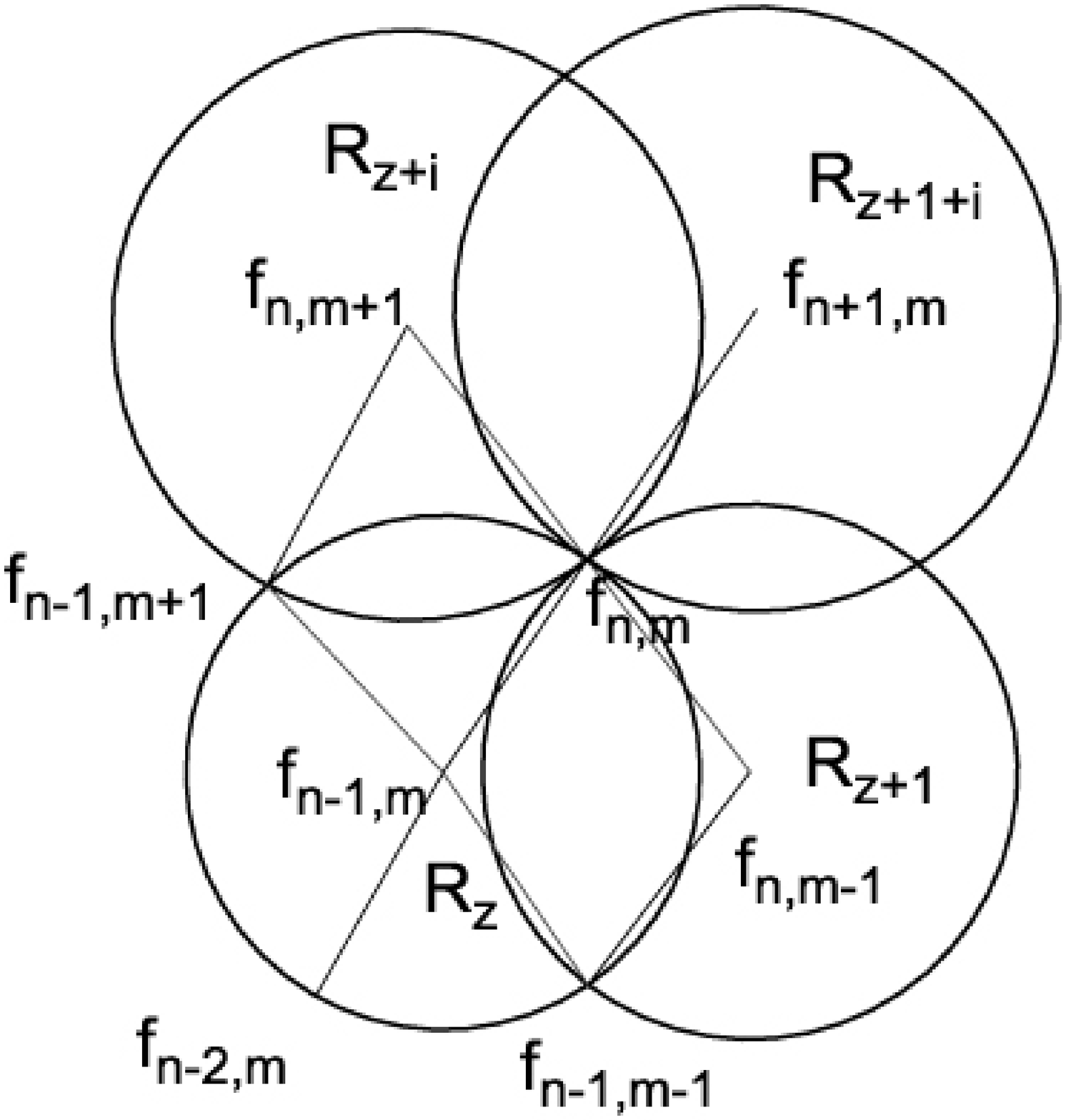,width=60mm} \ \ \ \ \ \ \ \
\epsfig{file=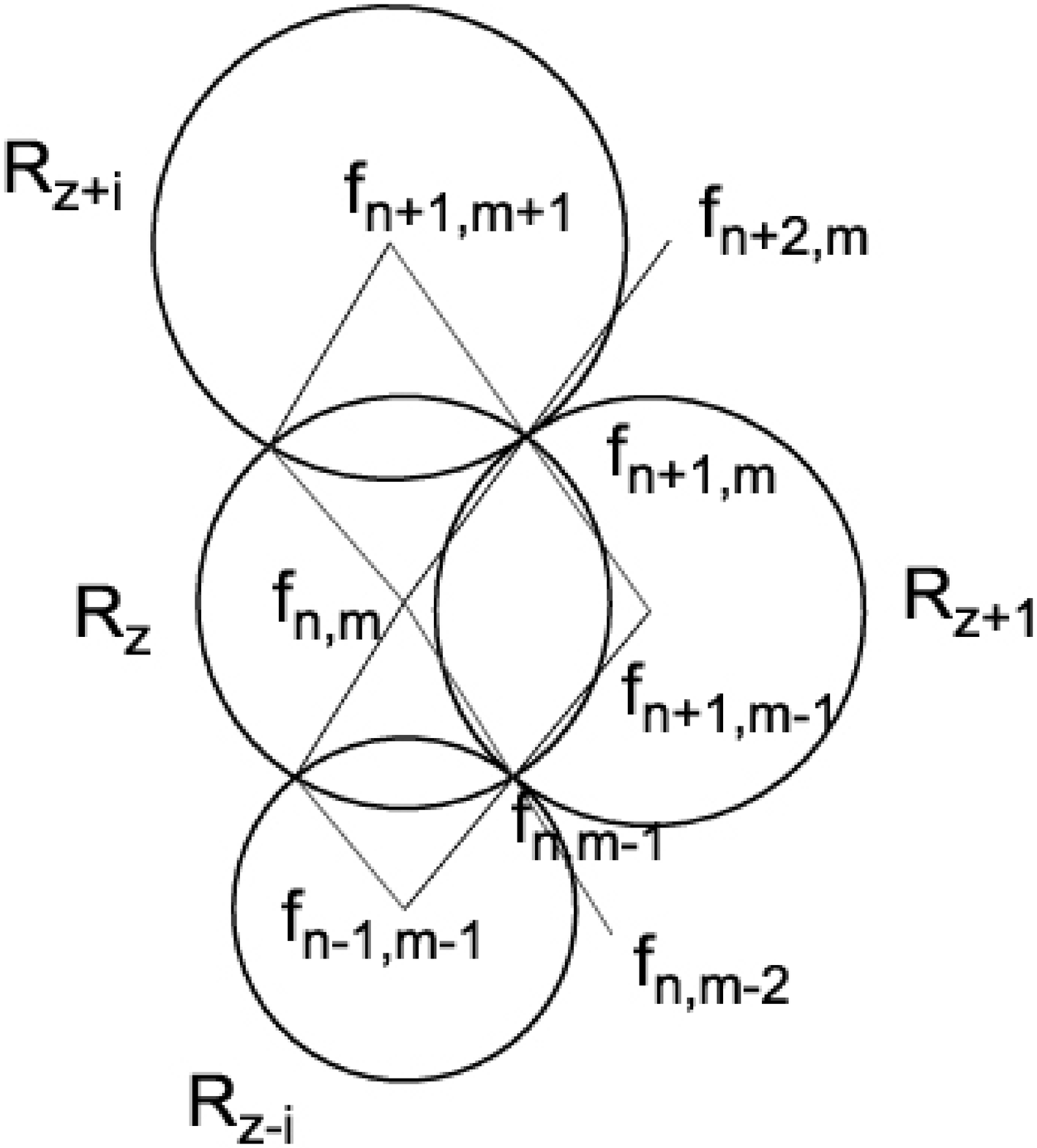,width=60mm} \caption{Kite-quadrilaterals of
circle pattern.
 } \label{REqSq}
\end{center}
\end{figure}

\begin{proposition}\label{eqforR}
Let the solution $f$ of (\ref{q}),(\ref{sg-c}) with initial
(\ref{skew-Shramm-initial}) is an immersion, then $R(z)$ defined
by (\ref{rmap}) satisfies the following equations:
\begin{equation}\label{square}
\small
 -MR_zR_{z+1}+(N+1)R_{z+1}R_{z+1+i}+(M+1)R_{z+1+i}R_{z+i}-NR_{z+i}R_z=\frac{\gamma}{2}(R_z+R_{z+1+i})(R_{z+1}+R_{z+i})
\end{equation}
for $z\in {\bf V}_l:={\bf V}\cup \{-N+i(N-1)|N\in {\bf N}\}$ and

\begin{equation}\label{Ri}
\small
(N+M)(R_{z+i}+R_{z+1})(R_z^2-R_{z+1}R_{z-i}+\cos \alpha
R_z(R_{z-i}-R_{z+1}))+
\end{equation}
$$ (M-N)(R_{z-i}+R_{z+1})(R_z^2-R_{z+1}R_{z+i}+\cos \alpha
R_z(R_{z+i}-R_{z+1}))=0, $$ for $z\in {\bf V}_{rint}:={\bf
V}\backslash \{\pm N+iN|N\in {\bf N}\}.$\\ Conversely let $R(z):
{\bf V} \to {\bf R_+}$ satisfy (\ref{square}) for $z\in {\bf V}_l$
and (\ref{Ri}) for $z\in {\bf V}_{rint}.$ Then $R(z)$ define a
square grid circle patterns with intersection angles $\alpha$, the
corresponding discrete  map $f_{n,m}$ is an immersion and
satisfies (\ref{q}),(\ref{sg-c}).
\end{proposition}
\noindent {\it Proof:} Circle pattern is immersed if and only if
all triangles $(f_{n,m},f_{n+1,m},f_{n,m+1})$ of  elementary
quadrilaterals of the map $f_{n,m}$ have the same orientation (for
brevity we call it orientation of quadrilaterals). Suppose that
the quadrilateral $(f_{0,0},f_{1,0},f_{1,1},f_{0,1})$ has positive
orientation. Let the circle pattern $f_{n,m}$ be an immersion. For
$n+m \equiv 1 \ ({\rm mod}\ 2) $  points
$f_{n,m},f_{n-1,m+1},f_{n-2,m},f_{n-1,m-1}$ lie on circle with the
center at $f_{n-1,m}$ and radius $R_z$, where
$z=(n-m-1)/2+i(n+m-1)/2$ (See the left part of Fig. \ref{REqSq}).
Using (\ref{q}) one can compute $f_{n,m+1}$ and $f_{n,m-1}$. Lemma
\ref{geom-kite} and proposition \ref{kite} imply that $f_{n+1,m}$
is in line with $f_{n-1,m},f_{n,m}$ and that the points
$f_{n,m+1},f_{n,m},f_{n,m-1}$ are also collinear. Denote by
$R_{z+1},R_{z+i}$ the radii of the circle at $f_{n,m-1}$ and
$f_{n,m+1}$ respectively. Define
$R_{z+1+i}=R_z\frac{(f_{n+1,m}-f_{n,m})}{(f_{n,m}-f_{n-1,m})}.$ If
(\ref{sg-c}) is satisfied at $(n-1,m)$ then (\ref{sg-c}) at
$(n,m)$ is equivalent to (\ref{square}), $R_{z+1+i}$ being
positive iff the quadrilaterals
$(f_{n,m},f_{n+1,m},f_{n+1,m+1},f_{n,m+1})$ and
$(f_{n,m-1},f_{n+1,m-1},f_{n+1,m},f_{n,m})$ have positive
orientation.

Similarly  starting with (\ref{sg-c}) at $(n,m-1)$, where $n+m
\equiv 0 \ ({\rm mod}\ 2) $ (see the right part of
Fig.\ref{REqSq}) one can determine evolution of the cross-like
figure formed by $f_{n,m-1}$, $f_{n+1,m-1}$, $f_{n,m}$,
$f_{n-1,m-1}$, $f_{n,m-2}$ into $f_{n+1,m}$, $f_{n+2,m}$,
$f_{n+1,m+1}$, $f_{n,m}$, $f_{n+1,m-1}$. Equation (\ref{sg-c}) at
$(n+1,m)$ is equivalent to (\ref{square}) and (\ref{Ri}) at
$z=(n-m)/2+i(n+m)/2$. $R_{z+1}$ is positive only for immersed
circle pattern.

Now let $R_z$ be some positive solution to
(\ref{square}),(\ref{Ri}). We can rescale it so that $R_0=1$. This
solution is completely defined by $R_0,R_i$. Consider solution
$f_{n,m}$ of (\ref{q}),(\ref{sg-c}) with initial data
(\ref{skew-Shramm-initial}) where $\beta$ is chosen so that the
quadrilateral $(f_{0,0},f_{1,0},f_{1,1},f_{0,1})$ has positive
orientation and satisfies the conditions $R_0=1=|f_{0,0}-f_{1,0}|$
and $R_i=|f_{1,1}-f_{1,0}|$. The map $f_{n,m}$ defines circle
pattern due to proposition \ref{kite}. It can be uniquely computed
from (\ref{square}),(\ref{Ri}). To this end one have to resolve
(\ref{square}) with respect to $R_{z+i+1}$ and use it to find
$f_{n+1,m}$ from
$R_{z+1+i}=\frac{R_z(f_{n+1,m}-f_{n,m})}{(f_{n,m}-f_{n-1,m}})$ and
to resolve (\ref{Ri}) for $R_{z+i}$ to find $f_{n+1,m+1}$ from
$R_{z+i}=R_{z+1}\frac{(f_{n+1,m+1}-f_{n+1,m})}{(f_{n+1,m}-f_{n+1,m-1})}$.
One can reverse the argument used in derivation of
(\ref{square}),(\ref{Ri}) to show that $f$ satisfies
(\ref{q}),(\ref{sg-c}). Moreover, since $R_{z}$ is positive, at
each step we get positively orientated quadrilaterals. Q.E.D.
\smallskip

Note that initial data (\ref{skew-Shramm-initial}) for $f_{n,m}$
imply initial data for $R_z$:
\begin{equation}\label{Rinitial}
R_0=1, \ \ \ R_i=\frac{\sin  \frac{\beta }{2}}{\sin (\alpha
-\frac{\beta}{2})}.
\end{equation}

\begin{definition}
A discrete  map $f_{n,m}$ is called embedded if inner parts of
different elementary quadrilaterals
$(f_{n,m},f_{n+1,m},f_{n+1,m+1},f_{n,m+1})$ do not intersect.
\end{definition}

\begin{proposition}\label{convex}
If for  a  solution $R_z$ of (\ref{square}),(\ref{Ri}) with
$\gamma \ne 1$ and initial conditions (\ref{Rinitial})  holds
true:
\begin{equation}
R_z>0, \ \ (\gamma -1)(R_z^2-R_{z+1}R_{z-i}+\cos \alpha
R_z(R_{z-i}-R_{z+1}))\ge 0 \label{sign}
\end{equation}
  in ${\bf V}_{int}$, then the
corresponding discrete  map is embedded.
\end{proposition}
\noindent {\it Proof:} Since $R(z)>0$ the corresponding discrete
map is an immersion.
  Consider piecewise linear curve $\Gamma _n$ formed by segments
$[f_{n,m},f_{n,m+1}]$, where $n>0$, $0\le m \le n-1$ and the
vector ${\bf v}_n(m)={ (f_{n,m}f_{n,m+1})}$ along this curve. Due
to Proposition \ref{kite}  this vector
 rotates only in vertices with
$n+m=0 \ ({\rm mod} \ 2)$ as $m$ increases along the curve. The
sign of the rotation angle $\theta _n(m)$, where $-\pi < \theta
_n(m) < \pi$, $0<m<n$ is defined by the sign of expression
$$R_z^2-R_{z+1}R_{z-i}+\cos \alpha R_z(R_{z-i}-R_{z+1})$$  (note
that there is no rotation if this expression vanishes), where
 $z=(n-m)/2+i(n+m)/2$ is a label for the circle with the center in
$f_{n,m}$. If $n+m=1\ ({\rm mod \ 2})$ define $\theta _n(m)=0$.

Now the theorem hypothesis and equation (\ref{Ri}) imply that
 the vector
${\bf v}_n(m)$ rotates with increasing $m$ in the same direction
for all $n$, and namely, clockwise for $\gamma <1$ and
counterclockwise for $\gamma >1$. Due to the following Lemma it is
sufficient to prove (\ref{sign}) only for $1<\gamma<2$.
\begin{lemma}\label{dual}
If $R_z$ is a solution of (\ref{square}),(\ref{Ri}) for $\gamma$
then $1/R_z$ is a solution of (\ref{square},\ref{Ri}) for $\tilde
\gamma =2-\gamma$. \label{dualR}
\end{lemma}

Consider the sector $B:=\{ z=re^{i\varphi}:r\ge 0,\ 0\le \varphi
\le \gamma \alpha /2 \}$. The terminal
 points of the curves $\Gamma _n$ lie
on the sector border.
\begin{lemma}\label{rot}
The curve $\Gamma _n$ has no self-intersection and lies in the
sector $B$.
\end{lemma}
Lemma is proved by induction. For $n=1$ it is obviously true since
the curve $\Gamma _1$ is a line segment. Let it be true for $n>1$.
Note that for immersed $f$ the first and the last segments of
$\Gamma _{n+1}$ lie in $B$ therefore the vector ${\bf v}_{n+1}(m)$
rotates counterclockwise around convex domain of $B$ confined by
$\Gamma _n$ and can not make a full circle. Lemma is proved.

\medskip

\noindent Each curve $\Gamma _n$ cuts the sector $B$ into a finite
part and an infinite part. Since the curve $\Gamma _n$ is convex
and the borders of all elementary quadrilaterals
$(f_{n,m},f_{n+1,m},f_{n+1,m+1},f_{n,m+1})$ for imbedded
$Z^{\gamma }$ have the positive orientation
 the segments of the
curve $\Gamma _{n+1}$ lie in the infinite part. Now the induction
in $n$ completes the proof of Proposition \ref{convex} for
$1<\gamma <2$. The proof for
 $0<\gamma <1$ is similar. Q.E.D.

\medskip

One can compute $Z^{\gamma }_{n,0}$ from (\ref{sg-c}) and obtain
the asymptotics $Z^{\gamma }_{n,0}\simeq c(\gamma)n^{\gamma}$ from
Stirling formula for large $x$:
\begin{equation}\label{Stirling}
\Gamma (x) \simeq \sqrt{2\pi }e^{-x}x^{x-\frac{1}{2}}.
\end{equation}
 Numerical experiments
support the following conjecture.

\medskip

\noindent {\bf Conjecture.} $Z^{\gamma }_{n,m}\simeq
c(\gamma)(n+e^{i\alpha}m)^{\gamma}$.

\section{Discrete Riccati equation and hypergeometric functions}\label{RE}
\begin{figure}[th]
\begin{center}
\epsfig{file=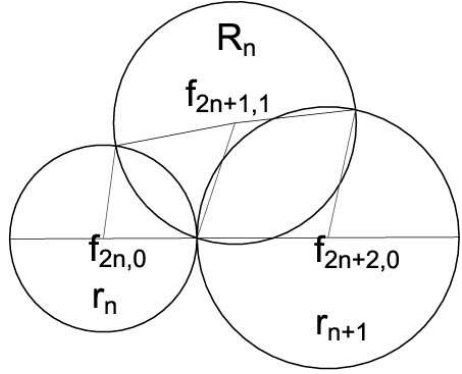,width=60mm} \caption{Circles on the border.
 } \label{RicC}
\end{center}
\end{figure}

Let $r_n$ and $R_{n}$ be radii of the circles with the centers at
$f_{2n,0}$, $f_{2n+1,1}$ respectively (see Fig.\ref{RicC}).

\noindent Constraint (\ref{sg-c}) and property (\ref{equidistant})
gives
\begin{equation}\label{diag-r}
r_{n+1}=\frac{2n+\gamma}{2(n+1)-\gamma}r_{n}.
\end{equation}
From elementary geometric considerations  one gets
$$R_{n+1}=\frac{r_{n+1}-R_n \cos \alpha}{R_n-r_{n+1}\cos
\alpha}r_{n+1}$$
 Define $$ p_n=\frac{R_{n}}{r_n}, \ \ \
g_n(\gamma)=\frac{2n+\gamma}{2(n+1)-\gamma}$$ and denote $t=\cos
\alpha$ for brevity. Now the equation for radii $R,\ r$ takes the
form:
\begin{equation}\label{Riccati}
p_{n+1}=\frac{g_n(\gamma )-tp_n}{p_n-tg_n(\gamma )}.
\end{equation}

\noindent {\bf Remark.} Equation (\ref{Riccati}) is a discrete
version of { \it  Riccati} equation. This title is motivated by
the following properties:
\begin{itemize}
\item cross-ratio of each four-tuple of its solutions is constant as $p_{n+1}$ is M\"obius transform of $p_n$,
\item general solution is expressed in terms of solution of
some linear equation (see below this linearisation).
\end{itemize}
 Below
we call (\ref{Riccati})  d-Riccati equation.

\begin{proposition} \label{R-positive}
Solution of discrete Riccati equation (\ref{Riccati})  is positive
for $n\ge 0$ iff
\begin{equation}\label{R-initial}
p_0=\frac{\sin  \frac{\gamma \alpha}{2}}{\sin \frac{(2-\gamma )
\alpha}{2}}
\end{equation}
\end{proposition}
Proof is based on the closed form of the general solution of
d-Riccati linearisation. It is linearised by the standard Ansatz
\begin{equation}\label{ansatz}
p_n=\frac{y_{n+1}}{y_n}+tg_n(\gamma )
\end{equation}
which transforms it into
\begin{equation}\label{linear}
y_{n+2}+t(g_{n+1}(\gamma )+1)y_{n+1}+(t^2-1)g_n(\gamma )y_n=0.
\end{equation}
One can guess that there is only one initial value $p_0$ giving
positive d-Riccati solution from the following consideration:
$g_n(\gamma )\to 1$ as $n\to \infty$, and the general solution of
 (\ref{linear}) with limit values of coefficients is
$y_n=c_1(-1)^n(1+t)^n+c_2(1-t)^n$. So
$p_n=\frac{y_{n+1}}{y_n}+tg_n(\gamma )\to -1$ for $c_1\ne 0$.
However  $c_1,c_2$ defines only asymptotics of a solution. To
relate it to initial values one needs some kind of {\it connection
formulas}. Fortunately it is possible to find the general solution
to (\ref{linear}).
\begin{proposition}\label{solution}
The general solution to (\ref{linear}) is
\begin{equation}\label{series}
y_n=\frac{\Gamma(n+\frac{1}{2})}{\Gamma (n+1-\frac{\gamma
}{2})}(c_1\lambda_1 ^{n+1-\gamma /2}F(\frac{3-\gamma
}{2},\frac{\gamma  -1}{2},\frac{1}{2}-n ,z_1)+
\end{equation}
$$ +c_2\lambda_2 ^{n+1-\gamma /2}F(\frac{3-\gamma
}{2},\frac{\gamma -1}{2},\frac{1}{2}-n ,z_2) $$ where $\lambda
_1=-t-1, \ \lambda _2=1-t, \ z_1 =(t-1)/2,\ z_2= -(1+t)/2$ and $F$
stands for the hypergeometric function.
\end{proposition}
\noindent {\it Proof:}  Solution was found by slightly modified
{\it symbolic method} (see \cite{Boole} for method description).
Substitution
\begin{equation}\label{substitution}
y_n=u_x\lambda ^x, \ \ x=n+1-\gamma /2
\end{equation}
transforms (\ref{linear}) into
\begin{equation}\label{x-linear}
\lambda ^2(x+1)xu_{x+2}+2t(x+\frac{\gamma
+1}{2})xu_{x+1}+(t^2-1)(x+\gamma -1)(x+1)u_x=0.
\end{equation}
We are looking for solution in the form
\begin{equation}\label{sum}
u_x=\sum_{m=-\infty}^{\infty}a_mv_{x,m}
\end{equation}
where $v_{x,m}$ satisfies
\begin{equation}\label{eq-for-v}
(x+m)v_{x,m}=v_{x,m+1}, \ \ xv_{x+1,m}=v_{x,m+1}.
\end{equation}

\noindent {\bf Remark.} Note that the label $m$ in (\ref{sum}) is
running by step 1 but is not necessary integer  therefore
$v_{x,m}$ is a straightforward generalization of
$x^{(m)}=(x+m-1)(x+m-2)...(x+1)x$ playing  the role of $x^m$ in
the calculus of finite differences. General solution to
(\ref{eq-for-v}) is expressed in terms of $\rm \Gamma$-function:
\begin{equation}\label{s-for-v}
v_{x,m}=c\frac{\Gamma(x+m)}{\Gamma(x)}
\end{equation}
Stirling formula (\ref{Stirling}) for large $x$ gives the
asymptotics for $v_{x,m}:$
\begin{equation}\label{asymptotics}
v_{x,m}\simeq cx^m \ \ {\rm for }\ \ x\to \infty .
\end{equation}
\medskip

Substituting (\ref{sum}) into (\ref{x-linear}), making use of
(\ref{eq-for-v}) and collecting similar terms one gets the
following equation for coefficients:
\begin{equation} \label{recurrent}
(\lambda ^2+2t\lambda+t^2-1)a_{m-2}+2(\frac{1+\gamma
}{2}-m)(t\lambda+t^2-1)a_{m-1}+(t^2-1)(1-m)(\gamma-1-m)a_m=0.
\end{equation}
Choice $\lambda _1=-t-1$ or $\lambda _2=1-t$ kills the term with
$a_{m-2}$.  To make series (\ref{sum}) convergent we can use the
freedom in $m$ to truncate (\ref{sum}) on one side. The choice
$m\in {\bf Z}$ or $m \in \gamma +{\bf Z}$ leads to divergent
series. For $m\in \frac{\gamma +1}{2}+{\bf Z}$ equation
(\ref{recurrent}) gives $a_{\frac{\gamma +1}{2}+k}=0$ for all
non-negative integer $k$ and
\begin{equation} \label{recurrent-1}
a_{\frac{\gamma
+1}{2}-k-1}=\frac{1-t^2}{t\lambda+t^2-1}\frac{(k-\frac{\gamma -1
}{2})(k-1+\frac{\gamma -1 }{2})}{2k}a_{\frac{\gamma +1}{2}-k}
\end{equation}
where $\lambda =\lambda _{1},\lambda _{2}$.   Substitution of
solution of this recurrent relation in terms of  $\Gamma
$-functions and (\ref{s-for-v}) yields
\begin{equation}\label{sum-1}
y_x=\lambda ^x
\sum_{k=1}^{\infty}\left(\frac{1-t^2}{2(t\lambda+t^2-1)}\right)^k
\frac{\Gamma (k-\frac{\gamma -1 }{2})\Gamma (k-1+\frac{\gamma -1
}{2})\Gamma (x+\frac{\gamma +1}{2}-k)}{\Gamma (k)\Gamma (x)}.
\end{equation}
\begin{lemma}\label{convergent}
For both $\lambda =-t-1,1-t$ series (\ref{sum-1}) converges for
all $x$.
\end{lemma}
{\it Proof of Lemma \ref{convergent}}: Since $z =
\frac{1-t^2}{2(t\lambda+t^2-1)}=(t-1)/2,-(1+t)/2$ for $\lambda
_{1},\lambda _{2}$ respectively and $t=\cos \alpha <1$ the
convergence of (\ref{sum-1}) depends on the behavior of
$\frac{\Gamma (k-\frac{\gamma -1 }{2})\Gamma (k-1+\frac{\gamma -1
}{2})\Gamma (x+\frac{\gamma +1}{2}-k)}{\Gamma (k)}$.  Stirling
formula (\ref{Stirling}) ensures  that this expression is bounded
by $ck^{\phi (x, \gamma )}$ for some $c$ an $\phi (x, \gamma )$
which gives convergence.
\smallskip

Series (\ref{sum-1}) is expressed in terms of hypergeometric
functions:
$$y_x=\lambda ^x\frac{\Gamma(x+\frac{\gamma
-1}{2})\Gamma (1-\frac{\gamma-1}{2})\Gamma
(\frac{\gamma-1}{2})}{\Gamma
(x)}F(1-\frac{\gamma-1}{2},\frac{\gamma-1}{2},1-\left(x+\frac{\gamma-1}{2}\right),z)
$$ where
\begin{equation} \label{hyper}
F(1-\frac{\gamma-1}{2},\frac{\gamma-1}{2},1-\left(x+\frac{\gamma-1}{2}\right),z)=1+z\frac{(1-\frac{\gamma-1}{2})(\frac{\gamma-1}{2})}{(1-(x+\frac{\gamma-1}{2}))}
+...+ \end{equation}
 $$
z^k\frac{
\left[(1-\frac{\gamma-1}{2})(2-\frac{\gamma-1}{2})...(k-\frac{\gamma-1}{2})\right]\left[(\frac{\gamma-1}{2})(1+\frac{\gamma-1}{2})...
(k-1\frac{\gamma-1}{2})\right]}{(1-(x\frac{\gamma-1}{2}))...(k-(x+\frac{\gamma-1}{2}))}+...
$$ Here we use the standard designation $F(a,b,c,z)$ for
hypergeometric function as a holomorphic at $z=0$ solution  for
equation
\begin{equation}\label{Gauss-eq}
z(1-z)F_{zz}+[c-(a+b+1)z]F_z-a b F=0.
\end{equation}

Now we can complete the proof of proposition \ref{solution}. Due
to linearity general solution of (\ref{linear}) is given by
superposition of any two linear independent solutions. As was
shown each summond  in (\ref{series})  satisfies the equation
(\ref{linear}). To finish the proof of Proposition \ref{solution}
one
 has to show that the particular solutions with $c_1=0,\ c_2\ne 0$ and $c_1 \ne 0,\ c_2=0$
 are linearly independent, which follows from Lemma \ref{as-lin}.
 \begin{lemma}\label{as-lin}
As $n\to \infty $ function (\ref{series})  has the asymtotics
\begin{equation}\label{eq-as-lin}
y_n \simeq (n+1-\gamma /2)^{\frac{\gamma -1}{2}}(c_1\lambda_1
^{n+1-\gamma /2}+c_2\lambda_2 ^{n+1-\gamma /2})
\end{equation}
 \end{lemma}
{\it Proof:} For $n\to \infty$ series representation (\ref{hyper})
gives $F(\frac{3-\gamma }{2},\frac{\gamma -1}{2},\frac{1}{2}-n
,z_1)\simeq1$. Stirling formula (\ref{Stirling}) defines
asymptotics of the factor $\frac{\Gamma(n+\frac{1}{2})}{\Gamma
(n+1-\frac{\gamma }{2})}$ (see (\ref{asymptotics}) ).Q.E.D.
\medskip

\noindent {\it Proof of Proposition \ref{R-positive}:} For
positive $p_n$ it is necessary that $c_1=0$: it follows from
asymptotics (\ref{eq-as-lin}) substituted into (\ref{ansatz}). Let
us define
\begin{equation}\label{s-def}
s(z)=1+z\frac{(1-\frac{\gamma -1}{2})(\frac{\gamma
-1}{2})}{\frac{1}{2}}+...+z^k\frac{(k-\frac{\gamma
-1}{2})...(1-\frac{\gamma -1}{2})(\frac{\gamma
-1}{2})(k-1+\frac{\gamma -1}{2})}{k!(k-\frac{1}{2})...\frac{1}{2}}
...
\end{equation}
It is the hypergeometric function $F(\frac{3-\gamma
}{2},\frac{\gamma -1}{2},\frac{1}{2}-n ,z)$ with $n=0$. A
straightforward manipulation with series shows that
\begin{equation}\label{p-0}
p_0=1+\frac{2(\gamma -1)}{2-\gamma }z+\frac{4z(z-1)}{2-\gamma
}\frac{s^\prime (z)}{s(z)}
\end{equation}
where $z=\frac{1+t}{2}$. Note that $p_0$ as a function of $z$
satisfies some ordinary differential equation of first order since
$\frac{s^\prime (z)}{s(z)}$ satisfies Riccati equation obtained by
reduction of (\ref{Gauss-eq}). Computation shows that $\frac{\sin
\frac{\gamma \alpha}{2}}{\sin \frac{(2-\gamma ) \alpha}{2}}$
satisfies the same ODE. Since both expression (\ref{p-0}) and
(\ref{R-initial}) are equal to 1 for $z=0$ they coincide
everywhere. Q.E.D.
\begin{corollary}\label{unique}
If there exists immersed  $f_{n,m}$ satisfying
(\ref{q}),(\ref{sg-c}),(\ref{skew-Shramm-initial}) it is defined
by initial data (\ref{right-sg-initial}).
\end{corollary}

\section{Embedded circle patterns and discrete Painlev\'e
equations}\label{IPE}

Let $R_z$ be a solution of (\ref{square}) and (\ref{Ri}) with
initial condition (\ref{Rinitial}). For $z\in {\bf V}_{int}$
define $P_{N,M}=P_z=\frac{R_{z+1}}{R_{z-i}}$,
$Q_{N,M}=Q_z=\frac{R_z}{R_{z-i}}$. Then (\ref{Ri}) and
(\ref{square})
 are rewritten as follows
\begin{equation}\label{QPainleve}
 Q_{N,M+1}= \frac{ (N-M)Q_{N,M}(1+P_{N,M})(Q_{N,M}-P_{N,M}\cos \alpha)-(M+N)P_{N,M}S_{N,M}
}{Q_{N,M}[ (M+N)S_{N,M}-(M-N)(1+P_{N,M})(P_{N,M}-Q_{N,M}\cos
\alpha ) ]}, \end{equation}
\begin{equation}
P_{N,M+1}=\frac{ (2M+\gamma )P_{N,M} + (2N+\gamma
)Q_{N,M}Q_{N,M+1}   }{ (2(N+1)-\gamma )P_{N,M} + (2(M+1)-\gamma
)Q_{N,M}Q_{N,M+1}  }, \label{PPainleve}
\end{equation}
where $$ S_{N,M}=Q_{N,M}^2-P_{N,M}+Q_{N,M}(1-P_{N,M})\cos \alpha
.$$ Property (\ref{sign}) for (\ref{QPainleve}),(\ref{PPainleve})
reads as
\begin{equation}
(\gamma-1)(Q_{N,M}^2 - P_{N,M}+Q_{N,M}(1-P_{N,M})\cos \alpha
)\ge 0, \ Q_{N,M}>0, \ P_{N,M}>0. \label{PQsign}
\end{equation}
Equations (\ref{QPainleve}),(\ref{PPainleve}) can be considered as
a dynamical  system  for  variable $M$, where due to
(\ref{diag-r})
\begin{equation}\label{initial-P}
P_{N+1,N}=\frac{2N+\gamma}{2(N+1)-\gamma}.
\end{equation}
\begin{proposition}\label{Q-exists}
For each $N\ge 0$ there exists $q_N>0$ such that the solution of
system (\ref{QPainleve}),(\ref{PPainleve}) subjected to
(\ref{initial-P}) and $Q_{N+1,N}=q_N$ has property (\ref{PQsign})
for all $M>N$.
\end{proposition}
\noindent {\it Proof:}   Due to Lemma \ref{dual} it is sufficient
to prove (\ref{PQsign}) only for $0<\gamma<1$. Define real
function $F(P)$ on ${\bf R_+}$ implicitly by $F^2 - P+F(1-P)\cos
\alpha = 0$ for $0\le P \le 1$ and by $F(P)\equiv 1$  for $1\le
P.$
\begin{figure}[th]
\begin{center}
\epsfig{file=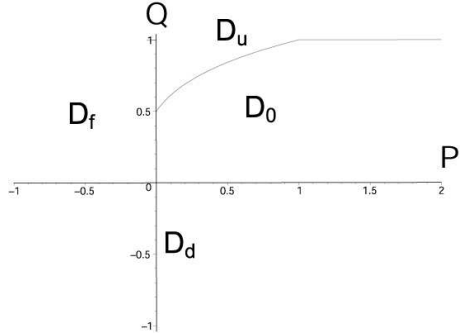,width=60mm} \caption{The case $\cos \alpha
=-1/2$.
 } \label{Domain}
\end{center}
\end{figure}

Designate
 $$
D_{u}:=\{(P,Q): P > 0, Q > F(P) \},\ \ D_{d}:=\{(P,Q): Q < 0 \},
$$ $$
 D_{0}:=\{(P,Q): P > 0, 0 \le
Q \le F(P) \}, \ \ D_f:=\{ (P,Q): P\le 0, Q\ge 0  \} $$ as in
Fig.\ref{Domain}.

System (\ref{QPainleve}),(\ref{PPainleve}) defines the map $\Phi_N
(M):(P_{N,M},Q_{N,M})\to(P_{N,M+1},Q_{N,M+1})$. This map has the
following properties:
\begin{itemize}
\item it is a continuous map on $D_{0}$. Values of $\Phi_N
(M)$ on the border of $D_{0}$ are defined by continuity in $\rm
\bf RP^2$.
\item For $(P_{N,M},Q_{N,M})\in D_{0}$ holds true $\Phi_N
(M,P_{N,M},Q_{N,M})\in D_{0}\cup D_{u} \cup D_{d}$, i.e. $(P,Q)$
can not jump in one step from $D_0$ into $D_f$.
\end{itemize}
Consider the solution to (\ref{QPainleve}),(\ref{PPainleve}) with
initial conditions on the segment $S(N)$ determined by
(\ref{initial-P}) and $Q_{N,N+1}=q$, where $0\le q \le
F(P_{N+1,N})$ then $(P_{N,M},Q_{N,M})=(P_{N,M}(q),Q_{N,M}(q))$.
Define $S_n(N)=\{q:(P_{N,M}(q),Q_{N,M}(q)) \in D_0 \forall  N < M
\le n\}$. Then $S_n(N)$ is a closed set as $\Phi_N (M)$ is
continuous on $D_0$ and $(P,Q)$ can not jump in one step from
$D_0$ on the half-line $P=0,Q>0$. As a closed subset of a segment
$S(N)$ the set $S_n(N)$ is a collection of disjoint segments
$\{S_n^l(N)\}$.
\begin{lemma} \label{segments}
There exists sequence $\{S_n^{l(n)}(N)\}_{n>N}$ such that:
\begin{itemize}
\item  $S_n^{l(n)}(N)$  is mapped by $(P_{N,n}(q),Q_{N,n}(q))$
onto some curve $c_n\in D_0$ with one terminal point on the curve
$Q=F(P)$ and the other on the line $Q=0$.
\item  $S_{n+1}^{l(n+1)}(N)\subset  S_n^{l(n)}(N)$.
\end{itemize}
\end{lemma}
Lemma is proved by induction. For $n=N+1$ it is trivial.  Let it
holds true for $n$. $\Phi_N (M)$ maps $(P,Q)\in D_0$ with $Q=F(P)$
into $D_u$ and $(P,0)$ with $P>0$ into $D_d$. Thus $\Phi_N (M)$
maps $c_n$ into some curve $\bar c_{n+1}$ lying in $D_{0}\cup
D_{u} \cup D_{d}$ with one terminal point  in $D_u$ and with the
other in $D_d$. Therefore at least one of the connected components
of $\bar c_{n+1}\cap D_0$ has its terminal points on the border of
$D_0$ as stated by Lemma.

\medskip

As the segments of $\{S_n^{l(n)}(N)\}$ constructed in Lemma
\ref{segments} are nonempty  there exits $q_N \in S_n(N)$ for all
$n> N$. For this $q_N$ holds true: $(P_{N,M}(q),Q_{N,M}(q))\in D_0
\setminus \partial D_0$. If $(P_{N,M}(q),Q_{N,M}(q))\in
\partial D_0$ then $(P_{N,M+1}(q),Q_{N,M+1}(q))$ would jump out of
$D_0$. Q.E.D.

\begin{theorem}\label{sg-imbed}
The discrete map $Z^{\gamma}$, $0<\gamma <2$ is embedded.
\end{theorem}
\noindent {\it Proof:}
 Proposition \ref{Q-exists} ensures that for
any $N>0$ there exist $R_{N+iN},R_{N+i(N+1)}$ such that the
solution to (\ref{square}),(\ref{Ri}) with these initial values is
positive for $z=N+iM$, $M\ge N$. Then equation (\ref{square})
defines recursively  $R_{K+iK}>0, \ K>N$ and implies that
$R_{n+im}>0$ far all $n,m: n\ge N,m\ge n$. Therefore asymptotics
(\ref{eq-as-lin}) implies that $R_{N+iN}/R_{N+i(N+1)}=Q_{N,N+1}$
is exactly as defined by initial condition (\ref{Rinitial}) for
$Z^{\gamma}$. Proposition \ref{convex} completes the proof.

\medskip

\noindent {\bf Remark 1. } For $N=0$ system
(\ref{QPainleve}),(\ref{PPainleve}) for $Q_{N,M},P_{N,M}$  reduces
 to the special case of discrete Painlev\'e equation (the case $\alpha =\pi/2$ was studied in \cite{A}):
 \begin{equation}
\label{dPII} (n+1)(x_n^2-1)\left(\frac{x_{n+1}+{x_n}/{\varepsilon}
}{\varepsilon+x_nx_{n+1}}\right)-
 n(1-{x_n^2}/{\varepsilon^2})\left(\frac{x_{n-1}+\varepsilon x_n}{\varepsilon+x_{n-1}x_n}\right) =
\gamma x_n \frac{\varepsilon ^2 -1}{2\varepsilon^2},
\end{equation}
where $\varepsilon =e^{i\alpha}$. This equation allows to
represent $x_{n+1}$ as a function of $n,x_{n-1}$ and $x_n$:
$x_{n+1}=\Phi (n,x_{n-1},x_n)$. $\Phi(n,u,v)$ maps the torus
$T^2=S^1\times S^1=\{(u,v)\in {\bf C}:\ |u|=|v|=1\}$ into $S^1$
and has the following properties:
\begin{itemize}
\item $\forall n\in {\bf N}$ it is a continuous map on $A_I\times
A_I$ where $A_I=\{e^{i\beta}: \beta \in [0,\alpha]\}$.
\item For $(u,v)\in A_I\times A_I$ holds true $\Phi(n,u,v)\in A_I\cup A_{II} \cup
A_{IV}$, where $A_{II}=\{e^{i\beta}: \beta \in (\alpha,\pi]\}$ and
$A_{IV}=\{e^{i\beta}: \beta \in [\alpha-\pi,0)\}$, i.e. $x$ can
not jump in one step from $A_I$ into $A_{III}=\{e^{i\beta}: \beta
\in (-\pi,\alpha-\pi )\}$.
\end{itemize}
These properties guarantees  that  there exists the unitary
solution $x_n=e^{i\alpha_n}$ of this equation with $x_0=e^{i\gamma
\alpha/2}$ in the sector $0<\alpha _n < \alpha$. This solution
corresponds to embedded $Z^{\gamma}$. Equation (\ref{dPII}) is a
special case of more general reduction of cross-ratio equation
(see \cite{N},\cite{AB} for the detail).

\section{Circle patterns $Z^2$ and $Log$}
For $\gamma=2$ formula (\ref{diag-r}) gives infinite $R_{1+i}$.
The way around this difficulty is re-normalization $f\to
(2-\gamma)f/\gamma$ and limit procedure $\gamma \to 2-0$, which
leads to the re-normalization of initial data (see \cite{BH}). As
follows from (\ref{R-initial}) this re-normalization gives:
\begin{equation} \label{2-initial}
R_0=0, \ R_{1+i}=1, \ R_i=\frac{\sin \alpha }{\alpha }.
\end{equation}
\begin{definition} Circle pattern $Z^{2}$ has radius function specified by
(\ref{square}),(\ref{Ri}) with initial data (\ref{2-initial}).
\end{definition}
 The
symmetry
\begin{equation} \label{duality} R \to \frac{1}{R},
\ \gamma \to 2-\gamma.
\end{equation} of (\ref{square}),(\ref{Ri}) stated in Lemma \ref{dual}
is the {\it duality transformation} (see \cite{BPD}). Smooth
analog $f\to f^*$ for holomorphic functions $f(w),f^*(w)$ is: $$
\frac{df(w)}{dw}\frac{df^*(w)}{dw}=1. $$ Note that ${\rm log }^*
(w)=w^2/2$.
\begin{definition}\cite{BH}
Circle pattern $\rm Log$ is a circle pattern dual to $Z^2$.
\end{definition}
\begin{theorem}
Discrete maps  corresponding to circle patterns $Z^2$ and $\rm
Log$ are embedded.
\end{theorem}
\noindent {\it Proof:} The circle radii for $Z^2$ and $\rm{Log}$
are subject to (\ref{square}),(\ref{Ri}) with $\gamma =2$ and
$\gamma=0$ respectively.  For these values of $\gamma$ Proposition
\ref{convex} is true: the proof is the same since $Z^2$ and
$\rm{Log}$ are immersed. Due to Lemma \ref{dual} it  suffices to
prove the property (\ref{sign}) only for $Z^2$.

Consider the discrete conformal map  $\frac{2-\gamma }{\gamma
}Z^{\gamma}$ with $0<\gamma <2.$  The corresponding
 solution $R_z$ of
 (\ref{square}),(\ref{Ri})
is a continuous function of $\gamma $. So there is   a limit  as
$\gamma \to 2-0$, of this solution with the property (\ref{sign}),
which is violated only for $z=0$ since $R_0=0$.\\ Q.E.D.

\section{Concluding remarks}\label{CR}
Further  generalizations of discrete $Z^{\gamma }$ and ${\rm Log}$
are possible.

\begin{figure}[th]
\begin{center}
\epsfig{file=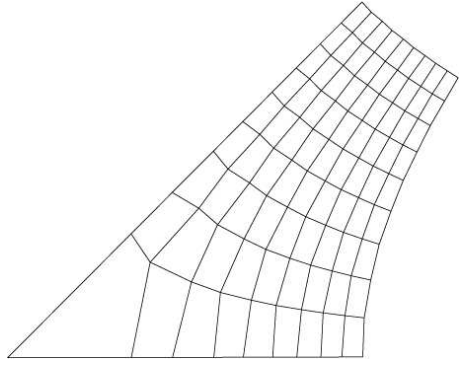,width=60mm} \caption{Discrete $Z^{1/2},\
\kappa =2,\ \alpha =\pi/2.$
 } \label{RecZ}
\end{center}
\end{figure}

\noindent {\bf I.} One can relax the unitary condition for
cross-ratios and consider solutions to
\begin{equation}\label{kq}
q(f_{n,m},f_{n+1,m},f_{n+1,m+1},f_{n,m+1})=\kappa ^2e^{-2i\alpha }
\end{equation}
subjected to the same constraint (\ref{sg-c}) with the initial
data
\begin{equation} \label{k-right-sg-initial}
  f_{1,0}=1,\
  f_{0,1}=\frac{e^{i \gamma \alpha  }}{\kappa}.
\end{equation}
This solution is a discrete analog of $Z^{\gamma}$ defined on the
vertices of regular parallelogram lattice (see Fig. \ref{RecZ}).
However, thus obtained mappings are deprived of geometrical flavor
as they do not define circle patterns.

\noindent {\bf II.} Another possibility is to de-regularize
prescribed combinatorics by projection of ${\bf Z^n}$ on a plane
as follows (see \cite{Se}). Consider ${\bf Z^n_+}\subset {\bf
R^n}$. For each coordinate vector ${\bf e_i}=(e^1_i,...,e^n_i)$
where $e^j_i=\delta ^j_i$ define unit vector ${\bf \xi _i}$ in
${\bf C}={\bf R^2}$ so that for any pair of indexes $i,j$ vectors
${\bf \xi_i},{\bf \xi_j}$ form a basis in ${\bf R^2}$. Let
$\Omega\in{\bf R^n}$ be some 2-dimensional connected simply
connected cell complex with vertices in ${\bf Z^n_+}$. Suppose
${0}\in \Omega$. (We denote by the same symbol $\Omega$ the set of
the complex vertices.)  Define the map $P:\Omega \to {\bf C}$ by
the following conditions:

1) $P(0)=0$,

2) if  $x,y$ are vertices of $\Omega$ and $y=x+{\bf e_i}$ then
$P(y)=P(x)+{\bf \xi _i}$.

\noindent It is easy to see that $P$ is correctly defined and
unique. \\ We call $\Omega$ a {\it projectable} cell complex iff
its image $\omega =P(\Omega)$ is embedded, i.e. intersections of
images of different cells of $\Omega$ do not have inner parts.
Using projectable cell complexes one can obtain not only regular
square grid and hexagonal combinatorics but more complex ones,
i.e. combinatorics of {\it Penrose tilings}.

It is natural to define discrete conformal map on $\omega$ as a
discrete complex immersion function $f$ on vertices of $\omega$
preserving cross-ratios of $\omega$-cells. The argument of $f$ can
be labeled by the vertices $x$ of $\Omega$. Hence for any cell of
$\Omega$, constructed on ${\bf e_k, e_j}$ the function $f$
satisfies the following equation for cross-ratios:
\begin{equation} \label{g-q}
q(f_x,f_{x+{\bf e_k}},f_{x+{\bf e_k}+{\bf e_j}},f_{x+{\bf
e_j}})=e^{-2i\alpha_{k,j} },
\end{equation}
where $\alpha_{k,j} $ is the angle between $\bf \xi_k$ and $\bf
\xi_j$, taken positive if $(\bf \xi_k, \bf \xi_j)$ has positive
orientation and taken negative otherwise. Now suppose that $f$ is
a solution to (\ref{g-q}) defined on the whole ${\bf Z^n_+}$.

\medskip

\noindent {\bf Conjecture 1.} {\it Equation (\ref{g-q}) is
compatible with the constraint }
\begin{equation} \label{g-c}
\gamma f_x=\sum_{s=1}^n2x_s \frac{(f_{x+{\bf
e_s}}-f_x)(f_x-f_{x-{\bf e_s}})}{f_{x+{\bf e_s}}-f_{x-{\bf e_s}}}
\end{equation}

\noindent For $n=3$ this conjecture was proven in \cite{BH}.

\noindent Now we can define discrete $Z^{\gamma}:\omega \to {\bf C
}$ for projectable $\Omega$ as solution to (\ref{g-q}),(\ref{g-c})
restricted on $\Omega$. Initial conditions for this solution are
of the form (\ref{right-sg-initial}) so that the restrictions of
$f$ on each two-dimensional coordinate lattices is an immersion
defining circle pattern with prescribed intersection angles.

This definition naturally generalize definition of discrete
$Z^{\gamma}$ given in \cite{BH} for $\Omega =\{(k,l,m):k+l+m=0,\pm
1 \}$.

\medskip

\noindent {\bf Conjecture 2.} {\it Discrete $Z^\gamma:\omega \to
{\bf C}$ is an immersion.}

\medskip

\noindent Schramm \cite{Schramm} showed that there is square grid
circle patterns mimicking ${\rm Erf} (\sqrt{i}z) $ but an analog
of ${\rm Erf} (z) $ does not exist. The obstacle is purely
combinatorial. There is a hope that combinatorics of projectable
cells can give more examples of discrete analogs of classical
functions.

\section{Acknowledgements}

The author thanks A.Bobenko, R.Halburd and A.Its for useful
discussions. This research was partially
 supported by the EPSRC grant
No Gr/N30941.

\end{document}